\documentclass{article}

\usepackage[T1]{fontenc}
\usepackage{geometry}			
\usepackage{parskip}			
\usepackage{mathtools}		
\usepackage[shortlabels]{enumitem}
\usepackage{amsthm}					
\usepackage{amssymb}	

\theoremstyle{plain}
	\newtheorem*{theorem*}{Theorem}
	\newtheorem{theorem}{Theorem}[section]
	\newtheorem{lemma}[theorem]{Lemma}
	
	\newtheorem{proposition}[theorem]{Proposition}

\theoremstyle{definition}
	
	\newtheorem{problem}{Problem}[section]

\theoremstyle{remark}
	
	\newtheorem*{notation}{Notation}

\newcommand{\card}[1]{\lvert #1 \rvert}
\newcommand{\floor}[1]{\lfloor #1 \rfloor}

\newcommand{\Set}[1]{\lbrace #1 \rbrace}
\newcommand{\defn}{\coloneq}			
\newcommand{\divides}{\mid}				
\newcommand{\transpose}[1]{#1^{T}}
\newcommand{\union}{\cup}

\newcommand{\intersect}{\cap}

\newcommand{\bigoh}[1]{O(#1)}					
\newcommand{\littleoh}[1]{o(#1)}				
\DeclareMathOperator{\rank}{rank}

\DeclareMathOperator{\charpol}{\mathsf{CharPol}}

\DeclareMathOperator{\Spec}{Spec}
\newcommand{\SymMat}{\mathsf{Sym}}

\usepackage{hyperref}	
	\hypersetup{
		colorlinks=true,
		linkcolor=red,
		urlcolor=blue,
		citecolor=magenta
	}
\usepackage[abbrev,nobysame,lite]{amsrefs}
\newcommand{\doi}[1]{\href{https://doi.org/#1}{doi:#1}}
\renewcommand{\MR}[1]{\href{https://mathscinet.ams.org/mathscinet-getitem?mr=#1}{MR#1}}
\newcommand{\Zbl}[1]{\href{https://zbmath.org/#1}{Zbl~#1}}
\newcommand{\arXiv}[1]{\href{https://arxiv.org/abs/#1}{arXiv:#1}}

\title{Low-rank matrices, tournaments, and symmetric designs\footnote{\copyright 2024. This manuscript version is made available under the CC BY-NC-ND 4.0 license \url{https://creativecommons.org/licenses/by-nc-nd/4.0/}. The published journal article is available at \emph{Linear Algebra Appl.}\ \textbf{694} (2024), 136--147, \href{https://doi.org/10.1016/j.laa.2024.04.006}{doi:10.1016/j.laa.2024.04.006}.}}

\author{
	Niranjan Balachandran
	\&
	Brahadeesh Sankarnarayanan\thanks{Corresponding author. Research supported by the National Board for Higher Mathematics (NBHM), Dept.\ of Atomic Energy (DAE), Govt.\ of India, and the Industrial Research and Consultancy Centre (IRCC), Indian Institute of Technology Bombay.}\\
	Department of Mathematics,
	Indian Institute of Technology Bombay,\\
	Mumbai 400076,
	Maharashtra,
	India.\\
	\texttt{\{niranj,bs\}@math.iitb.ac.in}
	}

\date{4 July 2026} 

\begin{document}
\maketitle
	
\begin{abstract}
Let \(\mathbf{a} = (a_{i})_{i \geq 1}\) be a sequence in a field \(\mathbb{F}\), and \(f \colon \mathbb{F} \times \mathbb{F} \to \mathbb{F}\) be a function such that \(f(a_{i},a_{i}) \neq 0\) for all \(i \geq 1\).
For any tournament \(T\) over \([n]\), consider the \(n \times n\) symmetric matrix \(M_{T}(f; \mathbf{a})\) with zero diagonal whose \((i,j)\)th entry (for \(i < j\)) is \(f(a_{i},a_{j})\) if \(i \to j\) in \(T\), and \(f(a_{j},a_{i})\) if \(j \to i\) in \(T\).
It is known~\cite{BalachandranBhattacharyaEtAl2023a} that if \(T\) is a uniformly random tournament over \([n]\), then \(\rank(M_{T}(f;\mathbf{a})) \geq (\frac{1}{2}-o(1))n\) with high probability when \(\operatorname{char}(\mathbb{F}) \neq 2\) and \(f\) is a linear function.

In this paper, we investigate the other extremal question: how low can the ranks of such matrices be?
We work with sequences \(\mathbf{a}\) that take only two distinct values, so the rank of any such \(n \times n\) matrix is at least \(n/2\).
First, we show that the rank of any such matrix depends on whether an associated bipartite graph has certain eigenvalues of high multiplicity.
Using this, we show that if \(f\) is linear, then there are \(n \times n\) real matrices \(M_{T}(f;\mathbf{a})\) of rank at most \(\frac{n}{2} + \bigoh{1}\).
For rational matrices, we show that for each \(\varepsilon > 0\) we can find a sequence \(\mathbf{a}(\varepsilon)\) for which there are \(n \times n\) matrices \(M_{T}(f;\mathbf{a}(\varepsilon))\) of rank at most \((\frac{1}{2} + \varepsilon)n + \bigoh{1}\).
These matrices are constructed from symmetric designs, and we also use them to produce bisection-closed families of size greater than \(\floor{3n/2} - 2\) for \(n \leq 15\), which improves the previously best known bound given in~\cite{BalachandranMathewEtAl2019}.
\end{abstract}

\section{Introduction}\label{S:Introduction}

Let \(a\) and \(b\) be two nonzero real numbers, and let \(T\) be a tournament over \([n]\).
Consider the \(n \times n\) symmetric matrix \(M_{T}\) with zero diagonal and \(M_{T}(i,j) = a\) if \(i \to j\) in \(T\), and \(M_{T}(i,j) = b\) if \(j \to i\) in \(T\), for all \(i < j\).

What can be said about the rank of \(M_{T}\)?

\subsection{Fractionally intersecting families: A motivating extremal problem}\label{SS:Motivation}

Linear algebraic methods are an important set of tools in the study of extremal problems in combinatorics; for more on this, see Babai--Frankl~\cite{BabaiFrankl2022}.
It is often the case that one can recast an extremal problem into the computation of certain matrix properties (its spectrum, rank, and so on).

The above question initially arose from the following extremal problem posed by Balachandran--Mathew--Mishra~\cite{BalachandranMathewEtAl2019}.
A \emph{\(\theta\)-intersecting family} \(\mathcal{F}\) over \([n]\) is a collection of subsets of \([n]\) such that for all distinct \(A,B \in \mathcal{F}\), we have \(\card{A \intersect B} = \theta \card{A}\) or \(\theta \card{B}\), where \(\theta \in (0,1) \intersect \mathbb{Q}\).
What is the maximum size of a \(\theta\)-intersecting family over \([n]\)?

The general problem remains open;
see~\cites{MathewMishraEtAl2022,BalachandranBhattacharyaEtAl2023c,Mishra2024} for related work.
A linear upper bound on the size of the \(\theta\)-intersecting family \(\mathcal{F} = \Set{A_{1},\dotsc,A_{m}}\) over \([n]\) can be obtained%
\footnote{The argument is implicit in~\cite{BalachandranMathewEtAl2019} and appears explicitly in~\cite{Sankarnarayanan2024}.}
by showing that an associated real symmetric \(m \times m\) matrix with zero diagonal, \(M_{\mathcal{F}}\), has ``high'' rank.%
\footnote{To be precise, if \(\rank(M_{\mathcal{F}}) \geq cm\), then \(\card{\mathcal{F}} \leq c^{-1}(n+1)\).}
Consequently, we turn our attention to the following ensemble of matrices.

\begin{notation}\label{N:ensemble}
	Let \(\mathbb{F}\) be a field, \(f \colon \mathbb{F} \times \mathbb{F} \to \mathbb{F}\) be a function, and \(a_{1},\dotsc,a_{n} \in \mathbb{F}\).
	By \(\SymMat(f;a_{1},\dotsc,a_{n})\) we denote the collection of all the \(n \times n\) symmetric matrices over \(\mathbb{F}\) with zero diagonal whose \((i,j)\)th entries are either \(f(a_{i},a_{j})\) or \(f(a_{j},a_{i})\), for all \(i < j\).
\end{notation}

\begin{notation}
Let \(\mathbb{F}\) be a field, \(\mathbf{a}\) be a sequence in \(\mathbb{F}\), and \(f \colon \mathbb{F} \times \mathbb{F} \to \mathbb{F}\) be a function.
\begin{enumerate}
	\item By \(\mathbf{a}_{n}\) we mean the tuple \((a_{1},\dotsc,a_{n})\).

	\item We say that \((f;\mathbf{a})\) is a \emph{good pair} (over \(\mathbb{F}\)) if \(f(a_{i},a_{i}) \neq 0\) for all \(i\).
	
	\item For \(\theta \in \mathbb{F}\), we write \(f = f_{\theta}\) if \(f(x,y) = x + (1-2\theta)y\) for all \(x,y \in \mathbb{F}\).
\end{enumerate}
\end{notation}

\begin{problem}\label{P:1}
	Let \((f;\mathbf{a})\) be a good pair over \(\mathbb{F}\).
	Is there an absolute constant \(c > 0\) such that \(\rank(M) \geq cn\) for all \(M \in \SymMat(f;\mathbf{a}_{n})\)?
\end{problem}

\subsection{\texorpdfstring{\(\SymMat(f;\mathbf{a}_{n})\)}{Sym(f; a n)}: A high-rank ensemble}\label{SS:Sym}

Problem~\ref{P:1} was initially posed in~\cite{BalachandranMathewEtAl2019} for \(\mathbb{F} = \mathbb{R}\), \(f = f_{1/2}\) and a sequence \(\mathbf{a}\) of positive reals.
In particular, the \(m \times m\) matrix \(M_{\mathcal{F}}\) mentioned earlier is a member of \(\SymMat(f_{\theta};\card{A_{1}},\dotsc,\card{A_{m}})\).%
\footnote{Specifically, it is the matrix whose \((i,j)\)th entry is \(f_{\theta}(\card{A_{i}},\card{A_{j}})\) if \(\card{A_{i} \intersect A_{j}} = \theta \card{A_{j}}\), and \(f_{\theta}(\card{A_{j}},\card{A_{i}})\) otherwise.}
However, this ensemble will in general contain many more matrices than just (those similar to) \(M_{\mathcal{F}}\), and we describe a natural combinatorial interpretation for all the matrices in \(\SymMat(f;\mathbf{a}_{n})\).

For a tournament \(T\) over \([n]\), define \(M_{T} \in \SymMat(f;\mathbf{a}_{n})\) by
\[
	M_{T}(i,j) =
		\begin{cases}
			f(a_{i},a_{j}), & i \to j \text{ in } T;\\
			f(a_{j},a_{i}), & j \to i \text{ in } T,
		\end{cases}
\]
for all \(i < j\).
Conversely, given a matrix \(M \in \SymMat(f;\mathbf{a}_{n})\), define a tournament \(T_{M}\) over \([n]\) by setting \(i \to j\) or \(j \to i\) in \(T_{M}\), for \(i < j\), depending on whether \(M(i,j) = f(a_{i},a_{j})\) or \(f(a_{j},a_{i})\), respectively.
Note that the correspondence \(M \to T_{M}\) is not one-one since we could have \(f(a_{i},a_{j}) = f(a_{j},a_{i})\) for some \(i < j\).

While there are good reasons to believe that \(M_{\mathcal{F}}\) has high rank,%
\footnote{Indeed, we do not even know of any \(\theta\)-intersecting family over \([n]\) of size at least \(2n\).}
is there evidence to suggest that all the matrices in the ensemble \(\SymMat(f;\mathbf{a}_{n})\) are of high rank?
In~\cite{TaoVu2017}, Tao and Vu study the rank of a random symmetric matrix \(M_{n} = ((\xi_{ij}))\), where \(\xi_{ij}\) are all jointly independent (for \(i < j\)) and also independent of \(\xi_{ii}\) (which are also independent) with the additional property that for all \(i < j\) and all real \(x\), \(\mathbb{P}(\xi_{ij} = x) \leq 1 - \mu\) for some fixed constant \(\mu\).
They show that the spectrum of \(M_{n}\) is simple with high probability (\emph{whp});
so, a random matrix in \(\SymMat(f;\mathbf{a}_{n})\) has rank at least \(n-1\) with high probability in the case when all the \(a_{i}\) are distinct reals.

To see why we require \((f;\mathbf{a})\) to be a good pair, consider the following matrix
\[
	M_{0} = \begin{bmatrix}
		0 & 1 & 4 & \cdots & (n-1)^{2} \\
		1 & 0 & 1 & \cdots & (n-2)^{2} \\
		4 & 1 & 0 & \cdots & (n-3)^{2} \\
		\vdots & \vdots & \vdots & \ddots & \vdots \\
		(n-1)^{2} & (n-2)^{2} & (n-3)^{2} & \cdots & 0
	\end{bmatrix},
\]
whose \((i,j)\)th entry is \((i-j)^{2}\).
It is easy to see that \(\rank(M_{0}) \leq 3\) for all \(n \geq 3\), since \(M_{0}\) is the sum of three rank one matrices, \(M_{1}\), \(M_{2}\), and \(M_{3}\), given by \(M_{1}(i,j) = i^{2}\), \(M_{2}(i,j) = -2ij\), \(M_{3}(i,j) = j^{2}\).
Note that \(M_{0} \in \SymMat(f;\mathbf{a}_{n})\), where \(f(x,y) = (x-y)^{2}\) and \(a_{i} = i\), but the same argument works for any sequence \(\mathbf{a}\) since \(f(x,x) = 0\) for all \(x \in \mathbb{F}\).
Hence, \(\rank(M) \leq 3\) for all \(M \in \SymMat(f;\mathbf{a}_{n})\) for \(f(x,y) = (x-y)^{2}\).

Another reason is that in~\cite{BalachandranBhattacharyaEtAl2023a} the authors showed that for a function \(f\) of finite rank, say \(f(x,y) = \sum_{i = 1}^{r} g_{i}(x)h_{i}(y)\) for some \(r \geq 1\), the rank of a matrix in \(\SymMat(f;\mathbf{a}_{n})\) arising from a uniformly random tournament is at least \((\frac{1}{2}-\littleoh{1})n\) \emph{whp}, provided that \((f;\mathbf{a})\) is a good pair over a field of characteristic different from \(2\).

\subsection{Low rank matrices in \texorpdfstring{\(\SymMat(f;\mathbf{a}_{n})\)}{Sym(f; a n)}}\label{SS:lowrank}

The discussion in the previous section shows why a typical member of \(\SymMat(f;\mathbf{a}_{n})\) has high rank.
However, we do not know of many concrete examples where the rank can be computed exactly.
In~\cite{BalachandranBhattacharyaEtAl2023b}, it is shown that the matrices in \(\SymMat(f_{1/2};\mathbf{a}_{n})\) associated to transitive tournaments all have high rank (specifically, at least \(n-1\)).
Furthermore, the probabilistic results mentioned in the previous section cannot rule out the existence of a vanishingly small proportion of matrices in \(\SymMat(f;\mathbf{a}_{n})\) whose ranks drop arbitrarily low.
So, it is interesting to see what can be said about matrices of \emph{low} rank in \(\SymMat(f;\mathbf{a}_{n})\).

The best examples that we know of are those associated to maximal bisection-closed families (i.e., \(\frac12\)-intersecting families).
The family \(\mathcal{F}_{\text{sunflower}}\)  over \([n]\) is constructed as a union of two sunflowers, containing \(n-1\) sets of size \(2\) and \(\frac{n}{2}-1\) sets of size \(4\).
The family \(\mathcal{F}_{\text{Hadamard}}\) over \([n]\) is constructed from a Hadamard matrix of order \(n\), and contains \(n-2\) sets of size \(n/2\) and \(\frac{n}{2}\) sets of size \(n/4\).
Since these maximal families over \([n]\) have size \(3n/2 - 2\), the associated matrices have rank at most \(\frac{2}{3}(3n/2 - 2) + O(1) = n + O(1)\).

Taking a cue from these constructions, we restrict our attention to matrices in \(\SymMat(f; \alpha^{(m)},\beta^{(n)})\), where \(\alpha,\beta \in \mathbb{R}\) satisfy \(f(\alpha,\alpha) \neq 0\) and \(f(\beta,\beta) \neq 0\), and \((\alpha^{(m)},\beta^{(n)})\) is shorthand for the sequence \((\alpha,\dotsc,\alpha,\beta,\dotsc,\beta)\) containing \(m\) occurrences of \(\alpha\) and \(n\) occurrences of \(\beta\).
Naturally, every matrix in \(\SymMat(f;\alpha^{(m)},\beta^{(n)})\) has rank at least \(r \defn \max\Set{m,n}\).
How close can one get to \(r\)?

Our first result says that the rank of a matrix \(M \in \SymMat(f;\alpha^{(m)},\beta^{(n)})\) depends on whether a scalar \(\mu(f;\alpha,\beta) \in \mathbb{C}\) such that
	\begin{equation}\label{Eq:mu}
		\bigl(\mu(f;\alpha,\beta)\bigr)^{2} = \frac{f(\alpha,\alpha)f(\beta,\beta)}{(f(\alpha,\beta)-f(\beta,\alpha))^{2}}
	\end{equation}
is an eigenvalue of high multiplicity for a bipartite graph \(G_{M}\) associated to \(M\), which we define in Section~\ref{S:Preliminaries}.

\begin{theorem}\label{T:Main}
	Let \(M \in \SymMat(f;\alpha^{(m)},\beta^{(n)})\), where \(f(\alpha,\alpha)\), \(f(\beta,\beta)\), and \(f(\alpha,\beta) - f(\beta,\alpha)\) are all nonzero.
	Let \(\mu(f;\alpha,\beta) \in \mathbb{C}\) be as in \eqref{Eq:mu}.
	There exists a bipartite graph \(G_{M}\) with part sizes \(m\) and \(n\) such that \(m + n - 2 - \nu \leq \rank(M) \leq m+n+2-\nu\), where \(\nu\) is the multiplicity of \(\mu(f;\alpha,\beta)\) as an eigenvalue of \(G_{M}\).
	In particular, if \(f(\alpha,\alpha) f(\beta,\beta) < 0\), then \(\rank(M) \geq m + n - 2\).
\end{theorem}

Note that when \(f = f_{\theta}\) (for \(\theta \in (0,1)\)), the hypotheses of Theorem~\ref{T:Main} are satisfied so long as \(\alpha\) and \(\beta\) are distinct and nonzero.
Next, we show that for any \(\theta \in (0,1)\) there exist \(\alpha,\beta \in \mathbb{R}\) for which there are matrices in \(\SymMat(f_{\theta};\alpha^{(n)},\beta^{(n)})\) having rank at most \(n+3\).

\begin{theorem}\label{T:2}
	Let \(\theta \in (0,1)\).
	Let \(\beta \in \mathbb{R}\) be a root of the quadratic equation \(x^{2} - (2 + (\theta^{-1} - 1)^{2})x + 1 = 0\).
	For every \(n \in \mathbb{N}\), there is a matrix \(M \in \SymMat(f_{\theta};1^{(n)},\beta^{(n)})\) with \(\rank(M) \leq n + 3\).
\end{theorem}

When \(\theta = 1/2\), the matrices that we find in the proof of Theorem~\ref{T:2} have irrational entries.
This raises the question whether there exist rational matrices in \(\SymMat(f_{1/2};\alpha^{(n)},\beta^{(n)})\) of rank close to \(n\).
While we are unable to settle the question completely, we show in the next theorem that we can get arbitrarily close to \(n\).

\begin{theorem}\label{T:Hadamard}
	For each \(\varepsilon > 0\), there exists \(c_{\varepsilon} \in (\frac12,\frac12+\varepsilon)\) and \(\alpha_{\varepsilon}, \beta_{\varepsilon} \in \mathbb{Q}\) such that there is a sequence of matrices \(M_{n} \in \SymMat(f_{1/2};\alpha_{\varepsilon}^{(n)},\beta_{\varepsilon}^{(n)})\) for which \(\rank(M_{n}) \leq 2c_{\varepsilon} n + \bigoh{1}\).
\end{theorem}

We conclude with some combinatorial consequences of these results.
In particular, we show that in some small cases (i.e., for \(n \leq 15\)), we can use these constructions to find bisection-closed families of size greater than \(\floor{3n/2} - 2\), which was the previously best known lower bound given in~\cite{BalachandranMathewEtAl2019}.

\section{Preliminaries}\label{S:Preliminaries}

For the rest of the paper, we will be working with real matrices in \(\SymMat(f;\alpha^{(m)},\beta^{(n)})\), where \(f(\alpha,\alpha) \neq 0\) and \(f(\beta,\beta) \neq 0\).
The latter conditions just say that \((f;\mathbf{a})\) is a good pair for any sequence \(\mathbf{a}\) taking values in \(\Set{\alpha,\beta}\).
Let \(I_{n}\) denote the \(n \times n\) identity matrix, let \(J_{m \times n}\) denote the \(m \times n\) all-ones matrix, and let \(J_{n} \defn J_{n \times n}\), for all \(m, n \in \mathbb{N}\).
Recall that by \(f_{\theta}\) we mean the function given by \((x,y) \mapsto x + (1-2\theta)y\).

To each matrix \(M \in \SymMat(f;\alpha^{(m)},\beta^{(n)})\), we associate a bipartite graph \(G_{M}\) with vertex set \(V(G_{M}) = \Set{v_{1},\dotsc,v_{m}} \union \Set{w_{1},\dotsc,w_{n}}\), as follows: for each \(i \in [m]\) and \(j \in [n]\), there is an edge \(\Set{v_{i},w_{j}}\) in \(G_{M}\) if and only if \(M(i,m+j) = f(\alpha,\beta)\).
Conversely, given a subgraph \(G\) of \(K_{m,n}\), the associated matrix \(M_{G} \in \SymMat(f;\alpha^{(m)},\beta^{(n)})\) is defined as follows:
if \(B_{m \times n}\) is the \(\Set{0,1}\)-incidence matrix of \(G\), let \(C = f(\alpha,\beta) B + f(\beta,\alpha) (J_{m \times n}-B)\), and define \(M_{G}\) to be the \((m+n) \times (m+n)\) matrix
\[
	\begin{bmatrix}
		f(\alpha,\alpha) (J_{m} - I_{m}) & C \\
		\transpose{C} & f(\beta,\beta) (J_{n} - I_{n})
	\end{bmatrix}.
\]
These correspondences are natural inverses of each other.

If \(f(\alpha,\beta) = f(\beta,\alpha)\), then there is only one matrix in the family \(\SymMat(f;\alpha^{(m)},\beta^{(n)})\), and it has rank at least \(m+n-2\).
So, to avoid this trivial scenario, we henceforth assume that \(f(\alpha,\beta) \neq f(\beta,\alpha)\).

The following lemma describes some properties of \(\mu(f_{\theta};\alpha,\beta)\), which we will refer to in Section~\ref{S:comcon}.
\begin{lemma}\label{L:sizes}
	If \(\alpha,\beta\) are distinct nonzero real numbers, then
	\begin{equation*}
		\bigl(\mu(f_{\theta};\alpha,\beta)\bigr)^{2} = \frac{(1-\theta)^{2}\alpha\beta}{\theta^{2}(\alpha-\beta)^{2}}.
	\end{equation*}
	In particular, \(\mu(f_{\theta};\alpha,\beta) = \mu(f_{\theta};\beta,\alpha)\) and \(\mu(f_{\theta};\alpha,\beta) = \mu(f_{\theta};k\alpha,k\beta)\) for all nonzero real numbers \(k\).
	Thus, \(\mu(f_{\theta};\alpha,\beta) = \mu(f_{\theta};1,\frac{\alpha}{\beta}) = \mu(f_{\theta};1,\frac{\beta}{\alpha})\).
\end{lemma}
\begin{proof}
	We have
	\begin{gather*}
		f_{\theta}(\alpha,\alpha) = 2(1-\theta)\alpha,\\
		f_{\theta}(\beta,\beta) = 2(1-\theta)\beta,\\
		f_{\theta}(\alpha,\beta)-f_{\theta}(\beta,\alpha) = 2\theta(\alpha-\beta).
	\end{gather*}
	Hence, 
	\[
		\bigl(\mu(f_{\theta};\alpha,\beta)\bigr)^{2} = \frac{f_{\theta}(\alpha,\alpha)f_{\theta}(\beta,\beta)}{(f_{\theta}(\alpha,\beta)-f_{\theta}(\beta,\alpha))^{2}} = \frac{(1-\theta)^{2}\alpha\beta}{\theta^{2}(\alpha-\beta)^{2}}.
	\]
	Since the RHS is symmetric and homogeneous of degree \(0\) in \(\alpha\) and \(\beta\), we have \(\mu(f_{\theta};\alpha,\beta) = \mu(f_{\theta};\beta,\alpha)\) and \(\mu(f_{\theta};\alpha,\beta) = \mu(f_{\theta};k\alpha,k\beta)\) for all nonzero real numbers \(k\).
	Taking \(k = 1/\alpha\), we get \(\mu(f_{\theta};\alpha,\beta) = \mu(f_{\theta};1,\frac{\alpha}{\beta})\).
	By symmetry, this equals \(\mu(f_{\theta};\frac{\alpha}{\beta},1)\).
	Now, taking \(k=\frac{\beta}{\alpha}\), we see that this also equals \(\mu(f_{\theta};1,\frac{\beta}{\alpha})\).
\end{proof}

We will refer to the following result due to Rowlinson in Section~\ref{SS:Rowlinson}, which concerns bipartite graphs having an eigenvalue of high multiplicity.
\begin{theorem}[Rowlinson~\cite{Rowlinson2016}*{Theorems 3.4 and 4.2}]\label{T:Rowlinson}
	Let \(G\) be a connected bipartite graph of order \(n > 5\), with \(\mu \notin \Set{-1,0}\) as an eigenvalue of multiplicity \(\nu > 1\).
	\begin{enumerate}[(a)]
		\item\label{(a)} If \(d\) is the maximum degree in \(G\), then \(\nu \leq n - 1 - d\).
		\item\label{(b)} If equality holds in \ref{(a)}, then \(\nu \leq d-1\).
		\item If equality holds in \ref{(b)}, then \(G\) is the bipolar cone over a graph \(G_{0}\), where \(G_{0}\) is either the incidence graph of a symmetric \(2\)-design, or a \(2\)-balanced bipartite graph.
	\end{enumerate}
\end{theorem}

\section{Main results}\label{S:Main}

\subsubsection*{Proof of Theorem~\ref{T:Main}}
	Without loss of generality, assume that \(m \geq n\).
	The matrix 
	\[
		\Gamma = \begin{bmatrix}
			f(\alpha,\alpha) J_{m} & f(\beta,\alpha) J_{m \times n}\\
			f(\beta,\alpha) J_{n \times m} & f(\beta,\beta) J_{n}
		\end{bmatrix}
	\]
	has rank at most \(2\), so \(\rank(M - \Gamma) - 2 \leq \rank(M) \leq \rank(M - \Gamma) + 2\).
	Since \(\rank(M - \Gamma) = m + n - \operatorname{nullity}(M - \Gamma)\), we shall find the multiplicity of \(0\) as a root of the characteristic polynomial \(\charpol(x)\) of \(M - \Gamma\):
	\[
		\charpol(x) = \det\begin{bmatrix}
			(x+f(\alpha,\alpha))I_{m} & (f(\beta,\alpha) - f(\alpha,\beta)) B\\
			(f(\beta,\alpha) - f(\alpha,\beta)) \transpose{B} & (x+f(\beta,\beta))I_{n}
		\end{bmatrix}.
	\]
	For a block matrix \(X = \begin{bmatrix} P & Q\\ R & S \end{bmatrix}\), if \(\det(P) \neq 0\), then \(\det(X) = \det(P) \det(S - RP^{-1}Q)\).
	So, \(\charpol(x)\) equals
	\[
		(f(\beta,\alpha) - f(\alpha,\beta))^{2n} (x+f(\alpha,\alpha))^{m-n} \det\left(\frac{(x+f(\alpha,\alpha))(x+f(\beta,\beta))}{(f(\beta,\alpha) - f(\alpha,\beta))^{2}}I_{n} - \transpose{B}B\right).
	\]
	Let \(A\) be the adjacency matrix of the bipartite graph \(G_{M}\), so \(A = \begin{bmatrix} 0 & B\\ \transpose{B} & 0 \end{bmatrix}\).
	Note that:
	\begin{itemize}
		\item the spectrum of \(A\) is real and symmetric about \(0\);
		\item \(A^{2} = \begin{bmatrix}B\transpose{B} & 0 \\ 0 & \transpose{B}B \end{bmatrix}\); and
		\item the spectra of \(B\transpose{B}\) and \(\transpose{B}B\) are identical (other than the \(m-n\) many additional zeros in the spectrum of \(B\transpose{B}\)).
	\end{itemize}
	Thus, \(0\) is an eigenvalue of \(M - \Gamma\) if and only if \(\mu(f;\alpha,\beta)\) is an eigenvalue of \(A\) (i.e., of \(G_{M}\)).\hfill$\square$

\subsubsection*{Proof of Theorem~\ref{T:2}}
	By Lemma~\ref{L:sizes}, there is no loss of generality in setting \(\alpha = 1\), so we state the result by specifying the required real number \(\beta\).
	It is easy to see that \(\Spec(K_{n}) = \Set{n-1, (-1)^{n-1}}\).
	Let \(G_{n}\) be the bipartite cover of \(K_{n}\), so that \(G_{n}\) is the complete bipartite graph \(K_{n,n}\) minus a perfect matching.
	Note that if \(v\) is an eigenvector of \(K_{n}\) with eigenvalue \(\lambda\), then \((v,v)\) and \((v,-v)\) are eigenvectors of \(G_{n}\) with eigenvalues \(\lambda\) and \(-\lambda\), respectively.
	Consequently, \(\Spec(G_{n}) = \Set{n-1, (1)^{n-1}, (-1)^{n-1}, -n+1}\).
	Given a nonzero \(\beta\), let \(M(\beta)\) be the matrix in \(\SymMat(f_{\theta};1^{(n)},\beta^{(n)})\) associated to \(G_{n}\).
	From Lemma~\ref{L:sizes}, we have
	\[
		\bigl(\mu(f_{\theta};1,\beta)\bigr)^{2} = \frac{(1-\theta)^{2}\beta}{\theta^{2}(1-\beta)^{2}}.
	\]
	So, we set \(\mu(f_{\theta};1,\beta) = 1\) and by rearranging terms we see that \(\beta\) must be a solution to the quadratic equation \(x^{2} - (2 + (\theta^{-1} - 1)^{2})x + 1 = 0\).
	In particular, we have
	\[
		\beta^{2} - (2 + (\theta^{-1} - 1)^{2})\beta + 1 = 0
		\iff \beta_{\pm} = \frac{2 + (\theta^{-1} - 1)^{2} \pm (\theta^{-1} - 1) \sqrt{(\theta^{-1} - 1)^{2} + 4}}{2}.
	\]
	Note that \(\beta_{\pm} \in \mathbb{R}\), so \(M(\beta_{\pm})\) is in fact a real matrix, and furthermore \(\rank(M(\beta_{\pm})) \leq n + n + 2 - (n-1) = n + 3\).\hfill$\square$

Note that \(\beta \in \mathbb{Q}\) for certain values of \(\theta\) in Theorem~\ref{T:2}; for example, if \(\theta = \frac25\), then \(\beta_{+} = 4\), and if \(\theta = \frac{3}{11}\), then \(\beta_{+} = 9\), and so on.
However, when \(\theta = \frac12\), we have \(\beta_{\pm} = \frac{3\pm\sqrt{5}}{2}\).
Nevertheless, can one find rational (or integral) matrices \(M \in \SymMat(f_{1/2};\alpha^{(n)},\beta^{(n)})\) having rank close to \(n\)?

To answer this question in the form of Theorem~\ref{T:Hadamard}, we turn to symmetric designs.
For \(\Delta\) a symmetric \(2\)-\((v,k,\lambda)\) design, let \(G_{\Delta}\) be its point-block bipartite incidence graph, and \(M_{\Delta} \in \SymMat(f;(\alpha^{(v)},\beta^{(v)}))\) be the matrix associated to \(G_{\Delta}\).%
\footnote{As an aside, the graph \(G_{M_{\Delta}}\) as defined in Section~\ref{S:Preliminaries} is precisely the graph \(G_{\Delta}\) in this setting.}
It is well-known that the spectrum of \(G_{\Delta}\) is given by \(\Set{ k, (\sqrt{k - \lambda})^{(v-1)}, (-\sqrt{k - \lambda})^{(v-1)}, -k }\).
So, Theorem~\ref{T:Main} implies that either \(\rank(M_{\Delta}) \geq 2v - 3\) or \(\rank(M_{\Delta}) \leq v + 3\).
We will use the matrices \(M_{\Delta}\) associated to Hadamard designs \(\Delta\) in order to prove Theorem~\ref{T:Hadamard}.
Recall that a Hadamard design is a symmetric \(2\)-\((v,k,\lambda)\) design with parameters \(v = 4n-1\), \(k = 2n-1\), and \(\lambda = n-1\), for \(n\) a positive integer.
A Hadamard design with these parameters exists if and only if there is a Hadamard matrix of order \(4n\).

\subsubsection*{Proof of Theorem~\ref{T:Hadamard}}
By Lemma~\ref{L:sizes}, there is no loss of generality in restricting to integral \(\alpha\) and \(\beta\).
The Paley construction~\cite{ColbournDinitz2007}*{\S V.6.1.6--7} gives a skew-Hadamard matrix of order \(q + 1\) for every prime power \(q \equiv 3 \pmod{4}\).
Also, Williamson~\cite{Williamson1944} shows that there is a Hadamard matrix of order \(\beta^{2} - \beta\) whenever there is a skew-Hadamard matrix of order \(\beta\).

So, given \(\varepsilon > 0\), fix \(N \in \mathbb{N}\) such that \(c_{\varepsilon} \defn \frac{4N+2}{8N-2} \in (\frac12,\frac12+\varepsilon)\).
Choose \(\beta_{\varepsilon}\) to be a prime power congruent to \(3\) modulo \(4\) that is sufficiently large so that \(n_{0} \defn \beta_{\varepsilon}^{2} - \beta_{\varepsilon} \geq N\), and let \(\alpha_{\varepsilon} \defn \beta_{\varepsilon} - 1\).
Let \(H\) be a normalized Hadamard matrix of order \(4n_{0}\) obtained from the Kronecker product of any two Hadamard matrices of orders \(n_{0}\) and \(4\), respectively.
Let \(\Delta\) be the Hadamard design induced by \(H\), and \(M_{\Delta} \in \SymMat(f_{1/2};\alpha_{\varepsilon}^{(4n_{0}-1)},\beta_{\varepsilon}^{(4n_{0}-1)})\) be the associated matrix.

By Theorem~\ref{T:Main}, \(\rank(M_{\Delta}) \leq 4n_{0} + 2\) if and only if \(\bigl(\mu(f_{1/2};\alpha_{\varepsilon},\beta_{\varepsilon})\bigr)^{2} = k - \lambda\), where \(k = 2n_{0}-1\) and \(\lambda = n_{0}-1\).
On the one hand, \(\bigl(\mu(f_{1/2};\alpha_{\varepsilon},\beta_{\varepsilon})\bigr)^{2} = (\beta_{\varepsilon} - 1)\beta_{\varepsilon}\).
On the other hand, by construction we have \((\beta_{\varepsilon} - 1)\beta_{\varepsilon} = n_{0}\), which equals \(k - \lambda\).
Thus, by Theorem~\ref{T:Main}, \(\rank(M_{\Delta}) \leq 4n_{0}+2\), and \(\frac{4n_{0}+2}{8n_{0}-2} \leq \frac{4N+2}{8N-2} = c_{\varepsilon}\).
Lastly, if \(G_{\Delta}\) is the bipartite graph associated to \(M_{\Delta}\), then we take arbitrarily many disjoint copies of \(G_{\Delta}\) to get a sequence of matrices \(M_{n} \in \SymMat(f_{1/2};(\beta_{\varepsilon} - 1)^{(n)},\beta_{\varepsilon}^{(n)})\) of rank at most \(2c_{\varepsilon} n + \bigoh{1}\), as required.\hfill$\square$

\section{Combinatorial consequences and concluding remarks}\label{S:comcon}

\subsection{Bipartite graphs with high multiplicity eigenvalues}\label{SS:Rowlinson}

In order to apply Theorem~\ref{T:Main} to find matrices of low rank in \(\SymMat(f;\alpha,\beta)\), we need to construct bipartite graphs that have some eigenvalue of high multiplicity (in particular, one equal to \(\mu(f;\alpha,\beta)\)).
The maximum multiplicity of any eigenvalue that is greater than \(-2\) for bipartite graphs (in fact, more generally for triangle-free graphs) was investigated by Rowlinson~\cite{Rowlinson2016}, and his results that we quoted in Section~\ref{S:Preliminaries} as Theorem~\ref{T:Rowlinson} show that if there is a bipartite graph having an eigenvalue of maximum possible multiplicity, then it arises as a bipolar cone either of the incidence graph of a symmetric design, or of a \(2\)-balanced bipartite graph.
Rowlinson further notes in~\cite{Rowlinson2016} that there are no known examples of the latter kind, so it is really the symmetric designs that are of interest when trying to construct bipartite graphs having an eigenvalue of high multiplicity.

\subsection{Coincidences in the set sizes of \texorpdfstring{\(\mathcal{F}_{\text{sunflower}}\)}{F sunflower} and \texorpdfstring{\(\mathcal{F}_{\text{Hadamard}}\)}{F Hadamard}}

As noted in Section~\ref{SS:lowrank}, the maximal bisection-closed families \(\mathcal{F}_{\text{sunflower}}\) and \(\mathcal{F}_{\text{Hadamard}}\) have sets of only two sizes (\(2\) and \(4\) in \(\mathcal{F}_{\text{sunflower}}\), and \(n/2\) and \(n/4\) in \(\mathcal{F}_{\text{Hadamard}}\)); moreover, the number of sets of size \(i\) in \(\mathcal{F}_{\text{sunflower}}\) and of size \(n/i\) in \(\mathcal{F}_{\text{Hadamard}}\), for \(i = 2,4\), are nearly equal.
Lemma~\ref{L:sizes} explains why this is more than just a coincidence.
We ignore small additive constants in this discussion for simplicity.

Suppose that \(\mathcal{F}_{1}\) and \(\mathcal{F}_{2}\) are two maximal \(\theta\)-intersecting families over \([N]\), such that each \(\mathcal{F}_{i}\) only contains sets of two distinct sizes, \(\alpha_{i}\) and \(\beta_{i}\).
Suppose that these families have the same size, so that the matrices \(M_{\mathcal{F}_{i}} \in \SymMat(f_{\theta};\alpha_{i}^{(m_{i})},\beta_{i}^{(n_{i})})\) have the same size.
As mentioned in Section~\ref{SS:Motivation}, these matrices will also have the same rank.
Theorem~\ref{T:Main} says that this is because \(\mu(f_{\theta};\alpha_{i},\beta_{i})\) is an eigenvalue of the same multiplicity for the associated graphs \(G_{i}\).

Now, one way for this to happen is that the graphs \(G_{i}\) are isomorphic to each other, and \(\mu(f_{\theta};\alpha_{1},\beta_{1}) = \mu(f_{\theta};\alpha_{2},\beta_{2})\).
As it turns out, the intersection patterns of the families \(\mathcal{F}_{\text{sunflower}}\) and \(\mathcal{F}_{\text{Hadamard}}\), which are encoded by the associated bipartite graphs, are \emph{almost identical}.%
\footnote{To be precise, there is a set \(A \in \mathcal{F}_{\text{Hadamard}}\) such that the graphs associated to the families \(\mathcal{F}_{1} = \mathcal{F}_{\text{Hadamard}} \setminus \Set{A}\) and \(\mathcal{F}_{2} = \mathcal{F}_{\text{sunflower}} \setminus \Set{12}\) are isomorphic.}
But what about the multiplicity of the eigenvalues?
Lemma~\ref{L:sizes} assures that \(\mu(f_{\theta};\alpha_{1},\beta_{1}) = \mu(f_{\theta};\alpha_{2},\beta_{2})\) if \(\alpha_{1}/\beta_{1} = (\alpha_{2}/\beta_{2})^{\pm1}\).
And, we do indeed have \(2/4 = (n/4)/(n/2)\) in comparing the sunflower and Hadamard families.

While the implications run in the opposite direction to the one we would ideally like, it puts the coincidences remarked upon above into a different light;
at the very least, it shows why there might exist these relationships between the sizes of the uniform subfamilies of \(\mathcal{F}_{\text{sunflower}}\) and \(\mathcal{F}_{\text{Hadamard}}\) as well as the sizes of the sets in them.

\subsection{Ruling out candidates for low rank matrices in \texorpdfstring{\(\SymMat(f_{1/2};\alpha^{(m)},\beta^{(n)})\)}{Sym(f 1/2; alpha m, beta n)}}

Theorem~\ref{T:Hadamard} shows that we can get a sequence \(M_{n}\) of integral matrices in \(\SymMat(f_{1/2};\alpha_{n}^{(n)},\beta_{n}^{(n)})\) whose ranks asymptotically tend to \(n\), using Hadamard designs.
We show that many of the other known families of symmetric designs cannot give us such matrices.
\begin{proposition}\label{Prop:onebytwo}
	Let \(\Delta\) be a symmetric \(2\)-\((v,k,\lambda)\) design with \(k - \lambda = p^{m}\) for some prime \(p\) and integer \(m \geq 1\).
	Let \(\alpha,\beta\) be distinct positive integers.
	Consider \(M_{\Delta} \in \SymMat(f_{1/2};\alpha^{(v)},\beta^{(v)})\).
	If \(\rank(M_{\Delta}) \leq v + 3\), then \(p^{m} = 2\).
\end{proposition}
\begin{proof}
	We may scale the matrix \(M_{\Delta}\) by \(1/\gcd(\alpha,\beta)\) without changing its rank, so we assume that \(\gcd(\alpha,\beta) = 1\).
	If \(\rank(M_{\Delta}) \leq v + 3\), then \(k - \lambda = \bigl(\mu(f_{1/2};\alpha,\beta)\bigr)^{2} = \frac{\alpha\beta}{(\alpha-\beta)^{2}}\).
	Now, it is easy to see that if there is a prime factor \(q\) of \((\alpha-\beta)^{2}\) that is also a factor of \(\alpha\beta\), then \(q\) must divide both \(\alpha\) and \(\beta\).
	Since \(\gcd(\alpha,\beta) = 1\), we must have \(\alpha = \beta \pm 1\), and \(p^{m} = k - \lambda = (\beta \pm 1)\beta\).
	Without loss of generality, let \(p^{m} = (\beta-1)\beta\).
	But, \(\gcd(\beta-1,\beta) = 1\), so \(\beta - 1 = 1\), that is, \(\beta = 2 = k - \lambda\).
	In particular, \(p^{m} = 2\).
\end{proof}

Proposition~\ref{Prop:onebytwo} rules out most of the known examples of symmetric designs as viable candidates for inducing matrices of low rank.
Indeed, a quick glance through~\cite{ColbournDinitz2007}*{\S II.6.8} shows that none of Families 1, 7, 8, 9, 13, or 14 are viable, except for the Fano plane \(PG(2,2)\) from Family 1, which we discuss in further detail below.
A similar analysis also rules out the Menon designs (Family 6 in~\cite{ColbournDinitz2007}*{\S II.6.8}) which have parameters \(v = 4t^{2}\), \(k = 2t^{2} - t\), \(\lambda = t^{2} - t\).

\subsection{\texorpdfstring{\(\theta\)}{theta}-intersecting families with sets of only two distinct sizes}

As mentioned in Section~\ref{SS:Motivation}, a positive answer to Problem~\ref{P:1} would prove that any \(\theta\)-intersecting family \(\mathcal{F}\) over \([n]\) has size at most linear in \(n\), whereas the best known upper bound~\cite{BalachandranMathewEtAl2019} is \(\card{\mathcal{F}} \leq O\bigl(\frac{n \log^{2} n}{\log \log n}\bigr)\).
When we restrict \(\mathcal{F}\) to have sets of only two distinct sizes, as in the known maximal families \(\mathcal{F}_{\text{sunflower}}\) and \(\mathcal{F}_{\text{Hadamard}}\), we have the trivial upper bound \(\card{\mathcal{F}} \leq 2n\), whereas \(\card{\mathcal{F}_{\text{sunflower}}} = \floor{3n/2} - 2 = \card{\mathcal{F}_{\text{Hadamard}}}\).
So, can the results in this work improve the trivial upper bound in this scenario?

Theorem~\ref{T:2} says we cannot hope to do so by considering ensembles \(\SymMat(f_{\theta};\alpha^{(m)},\beta^{(n)})\) of \emph{real} matrices.
Theorem~\ref{T:Hadamard} says that even restricting to rational or integral ensembles will not refine the trivial upper bound \(\card{\mathcal{F}} \leq 2n\).
Then, a further natural restriction would be to demand that the two set sizes be in \(1:2\) proportion, just as in \(\mathcal{F}_{\text{sunflower}}\) and \(\mathcal{F}_{\text{Hadamard}}\).
This means that we want to find matrices of low rank in \(\SymMat(f_{1/2};1^{(m)},2^{(n)})\).

\subsubsection*{The Fano plane}
The construction by Hadamard designs is no longer helpful here.
In fact, the following argument due to Gordon Royle~\cite{Royle2023} shows that the only symmetric \(2\)-\((v,k,\lambda)\) designs that give matrices of low rank in \(\SymMat(f_{1/2};1^{(n)},2^{(n)})\) are those arising from the Fano plane and its complement.
By Theorem~\ref{T:Main}, we require \(k - \lambda = \bigl(\mu(f_{1/2};1,2)\bigr)^{2} = 2\).
Any \(2\)-\((v,k,\lambda)\) design \(\Delta\) satisfies \((v-1)\lambda = k(k-1)\), so we must have \((v-1)\lambda = (\lambda + 2)(\lambda + 1)\).
That is, \(v = \lambda + 4 + \frac{2}{\lambda}\).
Since \(v,\lambda \in \mathbb{N}\), we must have \(\lambda \divides 2\).
\begin{enumerate}
	\item If \(\lambda = 1\), then \(v = 7\) and \(k = 3\).
	There is a unique design with these parameters, which is the Fano plane, and its incidence graph is the Heawood graph.
	\item If \(\lambda = 2\), then \(v = 7\), and \(k = 4\).
	Again, there exists a unique design with these parameters, which is the design whose blocks are the complements of the blocks of the Fano plane, and its incidence graph is the bipartite complement of the Heawood graph.
\end{enumerate}

The above constructions give matrices \(M_{n} \in \SymMat(f_{1/2};1^{(7n)},2^{(7n)})\), of order \(14n\), with rank at most \(\tfrac{4}{7} \cdot 14n + O(1) = 8n + O(1)\).
This improves the rank-to-order ratio from the \(2/3\) of the \(\mathcal{F}_{\text{sunflower}}\) and \(\mathcal{F}_{\text{Hadamard}}\) families to \(4/7\).
Furthermore, for \(n = 1\), consider the family \(\mathcal{F}_{\text{Fano}}\) over \([8]\):
\begin{align*}
	\mathcal{F}_{\text{Fano}} &= \mathcal{F}_{\text{sunflower}} \union \Set{1357,1368,1458,1467}\\
	&= \Set{12,13,14,15,16,17,18,1234,1256,1278,1357,1368,1458,1467}.
\end{align*}
\(\mathcal{F}_{\text{Fano}}\) is a bisection-closed family of size \(14\) over \([8]\), and it is the symmetric \(2\)-\((7,3,1)\) design, which induces the Heawood graph as its incidence graph \(G_{\text{Fano}}\).

Modifications to \(\mathcal{F}_{\text{sunflower}}\) as above using the sets \(1357\), \(1368\), \(1458\), and \(1467\) can produce bisection-closed families over \([n]\) of size more than \(\floor{3n/2} - 2\) up to \(n \leq 15\).
Also, it does not appear that there are any bisection-closed families arising from \(M_{n}\) for \(n > 1\).
As such, the original conjecture that any bisection-closed family over \([n]\) has size at most \(\floor{3n/2} - 2\) remains open for large \(n\).

\end{document}